\documentclass[11pt]{article}
\title{Minimal Surface Linear Combination Theorem}
\author{Michael Dorff and Stephen Taylor}

\usepackage{amsmath}
\usepackage{amssymb}
\usepackage{amsthm}
\usepackage{indentfirst}
\usepackage{latexsym}
\usepackage{graphicx}
\usepackage{graphics}
\usepackage{url}

\theoremstyle{plain}        \newtheorem{thm}{Theorem}
\theoremstyle{definition}   \newtheorem{defn}[thm]{Definiton}
\theoremstyle{definition}   \newtheorem{ex}{Example}
\theoremstyle{remark}       
\theoremstyle{plain}       \newtheorem{coro}{Corollary}
\theoremstyle{plain}        \newtheorem{lem}{Lemma}

\DeclareMathSymbol{\R}{\mathalpha}{AMSb}{"52}
\DeclareMathSymbol{\C}{\mathalpha}{AMSb}{"43}
\DeclareMathSymbol{\D}{\mathalpha}{AMSb}{"44}
\newcommand{\x}{\mathbf{x}}

\newcommand{\din}[2]{\displaystyle \int_{#1}^{#2}}
\newcommand{\limit}[2]{\displaystyle \lim_{{#1} \rightarrow {#2}} }

\def\re{\mathrm{Re}}
\def\im{\mathrm{Im}}

\def\bf{\textbf}

\begin{document}

\maketitle

\begin{abstract}
Given two univalent harmonic mappings $f_1$ and $f_2$ on $\mathbb{D}$,
which lift to minimal surfaces via the Weierstrass-Enneper
representation theorem, we give necessary and sufficient conditions
for $f_3=(1-s)f_1+sf_2$ to lift to a minimal surface for $s\in[0,1]$.  We then construct such mappings from
Enneper's surface to Scherk's singularly periodic surface, Sckerk's doubly periodic surface to the catenoid, and
the 4-Enneper surface to the 4-noid.
\end{abstract}

\section{Background}

Complex-valued harmonic mappings can be regarded as
generalizations of analytic functions. In particular, a harmonic
mapping is a complex-valued function $f=u+iv$, where the $C^2$
functions $u$ and $v$ satisfy Laplace's equation.  The Jacobian of
such a function is given by $J_f=u_xv_y-u_yv_x$. On a simply
connected domain $D\subset\mathbb{C}$, a harmonic mapping $f$ has
a canonical decomposition $f=h+\overline{g}$, where $h$ and $g$
are analytic in $D$, unique up to a constant \cite{css}. We will
only consider harmonic mappings that are univalent with positive
Jacobian on $\mathbb{D}=\{z:|z|<1\}$. The dilatation $\omega$ of a
harmonic map $f$ is defined by $\omega\equiv g'/h'$. A result by
Lewy \cite{lew} states that $|h'(z)|>|g'(z)|$ if and only if
$f=h+\bar g$ is sense--preserving and locally univalent. The
reader is referred to \cite{Duren} for many interesting results on
harmonic mappings.

One area of study is the construction of families of harmonic
mappings \cite{dms} and their corresponding minimal surfaces
\cite{bdp}, \cite{dorffszynal}, \cite{driverduren}. We now present some
necessary background concerning minimal surfaces. Let $M$ be an
orientable surface that arises from a differentiable mapping ${\bf
x}$ from a domain $V \subset \R ^2$ (or $\C$) into $\R ^3$, so
that ${\bf x}(u,v) =(x^1(u,v),x^2(u,v),x^3(u,v))$. The
parametrization $\x$ is isothermal (or conformal) if and only if
$\x_{u} \cdot \x_{v} = 0 \text{ and } \x_{u} \cdot \x_{u} = \x_{v}
\cdot \x_{v} (=\lambda>0)$. Note that there exists an isothermal
parametrization on any regular minimal surface (see \cite{dhkw}). Fix a
point $p$ on $M$. Let $\mathbf{t}$ denote a vector tangent to $M$
at $p$ and $\mathbf{n}$ the unit normal vector to $M$ at $p$. Then
$\mathbf{t}$ and $\mathbf{n}$ determine a plane that intersects
$M$ in a curve $\gamma$. The normal curvature $\kappa _\mathbf{t}$
at $p$ is defined to have the same magnitude as the curvature of
$\gamma$ at $p$ with the sign of $\kappa _\mathbf{t}$ chosen to be
consistent with the choice of orientation of $M$. The principal
curvatures, $\kappa _1$ and $\kappa _2$, of $M$ at $p$ are the
maximum and minimum of the normal curvatures $\kappa_\mathbf{t}$
as $\mathbf{t}$ ranges over all directions in the tangent space.
The {\it mean curvature} of $M$ at $p$ is the average value
$H=\frac{1}{2}(\kappa _1 + \kappa _2)$.

\begin{defn}
A {\it minimal surface} in $\R ^3$ is a regular surface for which
the mean curvature is zero at every point.
\end{defn}

The following standard theorem provides the link between harmonic
univalent mappings and minimal surfaces:

\begin{thm}(Weierstrass-Enneper Representation). Every regular minimal
surface has locally an isothermal parametric representation of the
form

\begin{align}
\label{WErep}
X & = \Big( \re \Big\{ \din{0}{z} p(1+q^2) dw \Big\}, \notag \\
  &  \hspace{0.75in} \Big. \re \Big\{
        \din{0}{z} -ip(1-q^2) dw \Big\}, \Big. \notag \\
  &  \hspace{1.2in} \Big.  \re \Big\{ \din{0}{z} -2ipq dw \Big\} \Big).
\end{align}
in some domain $D\subset\mathbb{C}$, where $p$ is analytic and $q$
is meromorphic in $D$, with $p$ vanishing only at the poles (if
any) of $q$ and having a zero of precise order $2m$ wherever $q$
has a pole of order $m$.  Conversely, each such pair of functions
$p$ and $q$ analytic and meromorphic, respectively, in a simply
connected domain $D$ generate through the formulas (\ref{WErep})
an isothermal parametric representation of a regular minimal
surface.
\end{thm}

We will use (\ref{WErep}) in the following form:

\begin{coro}
For a harmonic function $f=h+\overline{g}$, define the analytic
functions $h$ and $g$ by $h=\int^zpd\zeta$ and $g=\int^z
pq^2d\zeta$. Then the minimal surface representation
$(\ref{WErep})$ becomes

\begin{equation}
\label{wer2}
\left(\mathrm{Re}\{h+g\}, \; \;
\mathrm{Im}\{h-g\}, \; \;
2\mathrm{Im}\left\{\int_0^z\sqrt{h'g'}d\zeta\right\}\right)
\end{equation}
\end{coro}

\section{Harmonic Linear Combinations}

The main consideration of this work is the study of harmonic
mappings of the form $f_3=t f_1+ (1-t)f_2$, where $t \in [0,1]$
and $f_1$, $f_2$ are both harmonic mappings. We will provide
conditions for $f_3$ to lift to a minimal surface via (\ref{wer2}),
and demonstrate several examples which further the work of
\cite{dorffszynal} and relate seemingly disconnected minimal
surfaces. Let $f_1=h_1+\bar{g}_1$ and $f_2=h_2+\bar{g}_2$ be two
univalent harmonic mappings on $\mathbb{D}$, which lift to minimal
surfaces, with dilatations $q_1^2=g_1'/h_1'$ and $q_2^2=g_2'/h_2'$
respectively, where $q_1, q_2$ are analytic. Construct a third
harmonic mapping
\begin{align*}
f_3=& t f_1(z)+ (1-t) f_2 \\
    =& \big[t h_1(z)+ (1-t) h_2(z)\big]+\big[t \overline{g_1(z)}+
        (1-t) \overline{g_2(z)} \big] \\
    = & h_3+\overline{g_3}
\end{align*}
and define its dilatation to be $\omega_3=g'_3/h'_3$.

\begin{lem}
If $\omega_1=\omega_2$, then $\omega_3$ is a perfect square of an
analytic function and hence $f_3$ is locally univalent.
\end{lem}
\begin{proof}
Suppose that $\omega_1=\omega_2$.  Then we have
\begin{equation*}
\omega_3=\dfrac{t h_1'\omega_1+ (1-t) h_2'\omega_1}{t h_1'+ (1-t)
h_2'}=\omega_1,
\end{equation*}
which shows $\omega_3$ is a perfect square of an analytic
function. Since $f_1$ is univalent, $|\omega_3|=|q_1^2|>0$ and so
$f_3$ is locally univalent.
\end{proof}

We now seek to study conditions under which $f_3$ is globally
univalent and thus lifts to a minimal surface. To do his, we need
a few definitions and theorems.

\begin{defn}
A domain $D\subset\mathbb{C}$ is said to be convex in the
$e^{i\beta}$ direction if for all $a\in\mathbb{C}$ the set

$$D\cap\{a+te^{i\beta}:t\in\mathbb{R}\}$$

is either connected or empty.  Specifically, a domain is convex in
the direction of the imaginary axis if all lines parallel to the
imaginary axis have a connected intersection with the domain.
\end{defn}

\begin{thm}[\cite{convex}, \cite{convextwo}]
\label{thm:koepf}
Given a harmonic function $f=h+\overline{g}$,
let $\phi=h-g$. $\phi$ is convex in the $e^{i\beta}$ direction if

$$\mathrm{Re}\{\phi'(1+ze^{i(\alpha+\beta)})(1+ze^{-i(\alpha-\beta)})\}> 0$$

for some $\alpha\in\mathbb{R}$ and for all $z\in\mathbb{D}$ \label{unithm}
\end{thm}

The following theorem will allow us to prove global univalence of
a class of harmonic mappings.

\begin{thm}[Clunie and Sheil-Small, \cite{css}]
\label{convex} A harmonic function $f=h+\overline{g}$ locally
univalent in $U$ is a univalent mapping of $U$ onto a domain
convex in the $e^{i\beta}$ direction if and only if
$\phi=h-e^{i2\beta}g$ is a conformal univalent mapping of $U$ onto
a domain convex in the $e^{i\beta}$ direction.
\end{thm}

The following theorem  allows us to determine if a function maps
onto a domain convex in the direction of the imaginary axis:

\begin{thm}[Hengartner and Schober, \cite{hs}]
\label{hsthm} Suppose $f$ is analytic and non-constant in $\D$.
Then
\begin{equation*}
\re \{(1-z^2)f'(z)\} \geq 0, z \in \D
\end{equation*}
if and only if
\begin{enumerate}
\item $f$ is univalent in $\D$,
\item $f$ is convex in the imaginary direction, and
\item there exists points $z_n'$, $z_n''$ converging to $z=1$,
$z=-1$, respectively, such that
 \begin{equation}
 \label{primeends}
  \begin{split}
 & \limit{n}{\infty} \re \{f(z_n')\}=\sup_{|z|<1}  \re \{f(z)\} \\
  & \limit{n}{\infty} \re \{f(z_n')\}=\sup_{|z|<1} \re \{f(z)\}.
\end{split}
\end{equation}
\end{enumerate}
\end{thm}

Note that the the normalization in (\ref{primeends}) can be
thought of in some sense as if $f(1)$ and $f(-1)$ are the right
and left extremes in the image domain in the extended complex
plane.

Using the above results, we derive the following two theorems.

\begin{thm}
\label{thm:dorff}
Let $f_1=h_1+\overline{g_1}$,
$f_2=h_2+\overline{g_2}$ be harmonic mappings convex in the
imaginary direction. Suppose $\omega_1=\omega_2$ and
$\phi_i=h_i-g_i$ is univalent, convex in the imaginary direction,
and satisfies the normalization given in (\ref{primeends}) for
$i=1,2$. Then $f_3=tf_1+(1-t)f_2$ is convex in the imaginary
direction $(0 \leq t \leq 1)$.
\end{thm}
\begin{proof}
We want to show that $\phi_3=t\phi_1+(1-t)\phi_2$ is convex in the
imaginary direction. Then by Theorem \ref{convex}, $f_3$ is convex
in the imaginary direction. By the hypotheses, Theorem \ref{hsthm}
applies to $\phi_1$, $\phi_2$. That is,
\begin{equation*}
\re \{ (1-z^2) \phi_i'(z) \} \geq 0, \forall i=1,2.
\end{equation*}
Consider
\begin{align*}
\re \{(1-z^2)\phi'_3(z)\} =& \re \{ (1-z^2)
(t\phi_1'(z)+(1-t)\phi_2'(z))\} \\
=&t\re \{(1-z^2) \phi_1'(z)\} +(1-t)\re \{(1-z^2)\phi_2'(z) \}
\geq 0.
\end{align*}
Hence, by applying Theorem \ref{hsthm} again, $\phi_3$ is convex
in the imaginary direction.
\end{proof}

We need not only restrict to surfaces convex in the imaginary direction.  The following gives a condition for a
function to be convex in an arbitrary direction:

\begin{thm}
\label{thm:taylor} For a harmonic function $f=h+\overline{g}$,
define $h-g=\phi=\phi_R+i\phi_I$. Then $\phi$ is convex in the
$e^{i\beta}$ direction if

\begin{equation}
\left[\cos\alpha+\cos(\beta+\gamma)\right]\left[\phi_R'\cos(\beta+\gamma)
-\phi_I'\sin(\beta+\gamma)\right]>0
\end{equation}

for some $\alpha \in \R$ and for all $z=r
e^{i\gamma}\in\mathbb{D}$.
\end{thm}

\begin{proof}
This theorem follows by applying Theorem \ref{thm:koepf} to $\phi$
to get
\begin{align*}
& \mathrm{Re}\{(\phi'_R+i\phi_I')(1+r
e^{i(\alpha+\beta+\gamma)})(1+r e^{i(\gamma-\alpha+\beta)})\} \\
& \hspace{.5in} = \phi_R'+2 r
\cos\alpha(\phi_R'\cos\theta-\phi_I'\sin\theta)+r^2(\phi_R'\cos
2\theta-\phi_I'\sin 2\theta \\
& \hspace{.5in}  = 2
(\cos\alpha+\cos(\beta+\gamma))(\phi_R'\cos(\beta+\gamma)-\phi_I'\sin(\beta+\gamma))
> 0,
\end{align*}
where $\theta=\beta + \gamma$.
\end{proof}

\section{Examples}

We now proceed to give two interesting examples resulting from
Theorems \ref{thm:dorff} and \ref{thm:taylor}.

\begin{ex}[Ennepers to Scherks singly-periodic] \mbox{}

\noindent Consider the harmonic maps
\begin{align*}
f_E = & z+\dfrac{1}{3}\overline{z}^3 \\
f_S = & \bigg[ \dfrac{1}{4} \ln \bigg( \dfrac{1+z}{1-z} \bigg) +
\dfrac{i}{4} \ln \bigg( \dfrac{i-z}{i+z} \bigg) \bigg] \\
  & \hspace{0.8in} +
\overline{\bigg[ \dfrac{1}{4} \ln \bigg( \dfrac{1+z}{1-z} \bigg) -
\dfrac{i}{4} \ln \bigg( \dfrac{i-z}{i+z} \bigg) \bigg]}
\end{align*}

It is straight forward to show that their dilatations are $\omega
= z^2$ and both harmonic maps satisfy the hypotheses of Theorem
\ref{thm:dorff}. Hence
\begin{equation*}
f_t=(1-t)f_E+tf_S
\end{equation*}
is globally univalent on $z \in \D$ and $\forall t \in [0,1]$. By
Corollary \ref{wer2}, $f_t$ lifts to a family of minimal surfaces. Note
that $f_0$ lifts to Ennepers surface parametrized by:
\begin{equation*}
X_0= \Bigg( \re \bigg\{ z+\dfrac{1}{3}z^3 \bigg\}, \im \bigg\{
z-\dfrac{1}{3} z^3 \bigg\}, \im \bigg\{ z^2 \bigg\} \Bigg)
\end{equation*}
and $f_1$ lifts to Scherks singly-periodic surface parametrized by
\begin{equation*}
X_1 = \Bigg( \re \bigg\{ \dfrac{1}{2} \ln \bigg( \dfrac{1+z}{1-z}
\bigg) \bigg\}, \im \bigg\{ \dfrac{i}{2} \ln \bigg(
\dfrac{i-z}{i+z} \bigg) \bigg\}, \im \bigg\{ \dfrac{1}{2} \ln
\bigg( \dfrac{1+z^2}{1-z^2} \bigg) \bigg\} \Bigg).
\end{equation*}

So for $t \in [0,1]$ we get a continuous family of minimal
surfaces transforming from Ennepers to Scherks singly-periodic. In
Figure \ref{enntosch}, we have shown six equal increments in
this transformation.

\begin{center}
\begin{figure}[h]
\centerline{\hbox{
\includegraphics[height=5in]{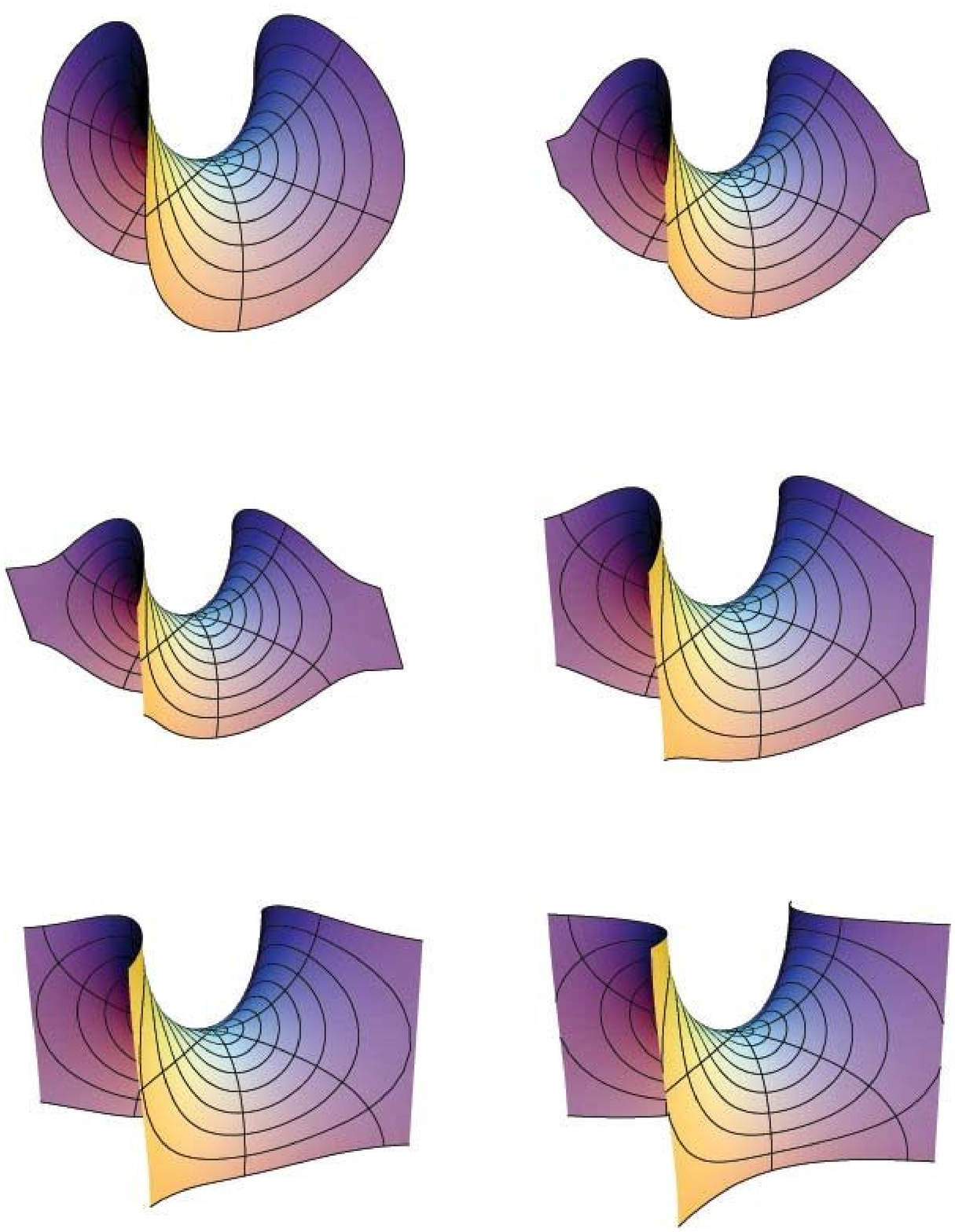}}}
\caption{Ennepers to Scherks singly-periodic transformation for
$t=i/5$ for $i=0,\ldots 5$.} \label{enntosch}
\end{figure}
\end{center}
\end{ex}

\begin{ex}[Scherks doubly-periodic to catenoid] \mbox{}

\noindent Consider the harmonic maps $f_D=h_D+\overline{g_D}$,
where
\begin{align*}
h_D(z)= &\dfrac{1}{4} \ln \bigg( \dfrac{1+z}{1-z} \bigg)
-\dfrac{i}{4} \ln \bigg( \dfrac{1+iz}{1-iz} \bigg) \\
g_D(z)=&-\dfrac{1}{4} \ln \bigg( \dfrac{1+z}{1-z} \bigg)
-\dfrac{i}{4} \ln \bigg( \dfrac{1+iz}{1-iz} \bigg),
\end{align*}
and $f_C=h_C+\overline{g_C},$ where
\begin{align*}
h_C(z)=\dfrac{1}{4} \ln \bigg( \dfrac{1+z}{1-z} \bigg)
+\dfrac{1}{2}
\dfrac{z}{1-z^2} \\
g_C(z)=\dfrac{1}{4} \ln \bigg( \dfrac{1+z}{1-z} \bigg)
-\dfrac{1}{2} \dfrac{z}{1-z^2}.
\end{align*}
Notice that both $f_D$ and $f_C$ are convex in the direction of
the imaginary axis, satisfy the hypotheses of Theorem
\ref{thm:dorff}, and for each $\omega=-z^2$. Hence
\begin{equation*}
f_t=(1-t)f_D+tf_C
\end{equation*}
is globally univalent on $z \in \D$ and $\forall t \in [0,1]$.

Note that $f_t$ lifts to a family of minimal surfaces, where $f_0$ lifts to
Scherks doubly-periodic surface and $f_1$ lifts to a catenoid. So
for $t \in [0,1]$ we get a continuous family of minimal surfaces
transforming from Scherks doubly-periodic surface to a catenoid.
In Figure \ref{schtocat}, we have shown six equal increments
in this transformation.

\begin{center}
\begin{figure}[h]
\centerline{\hbox{
\includegraphics[height=5in]{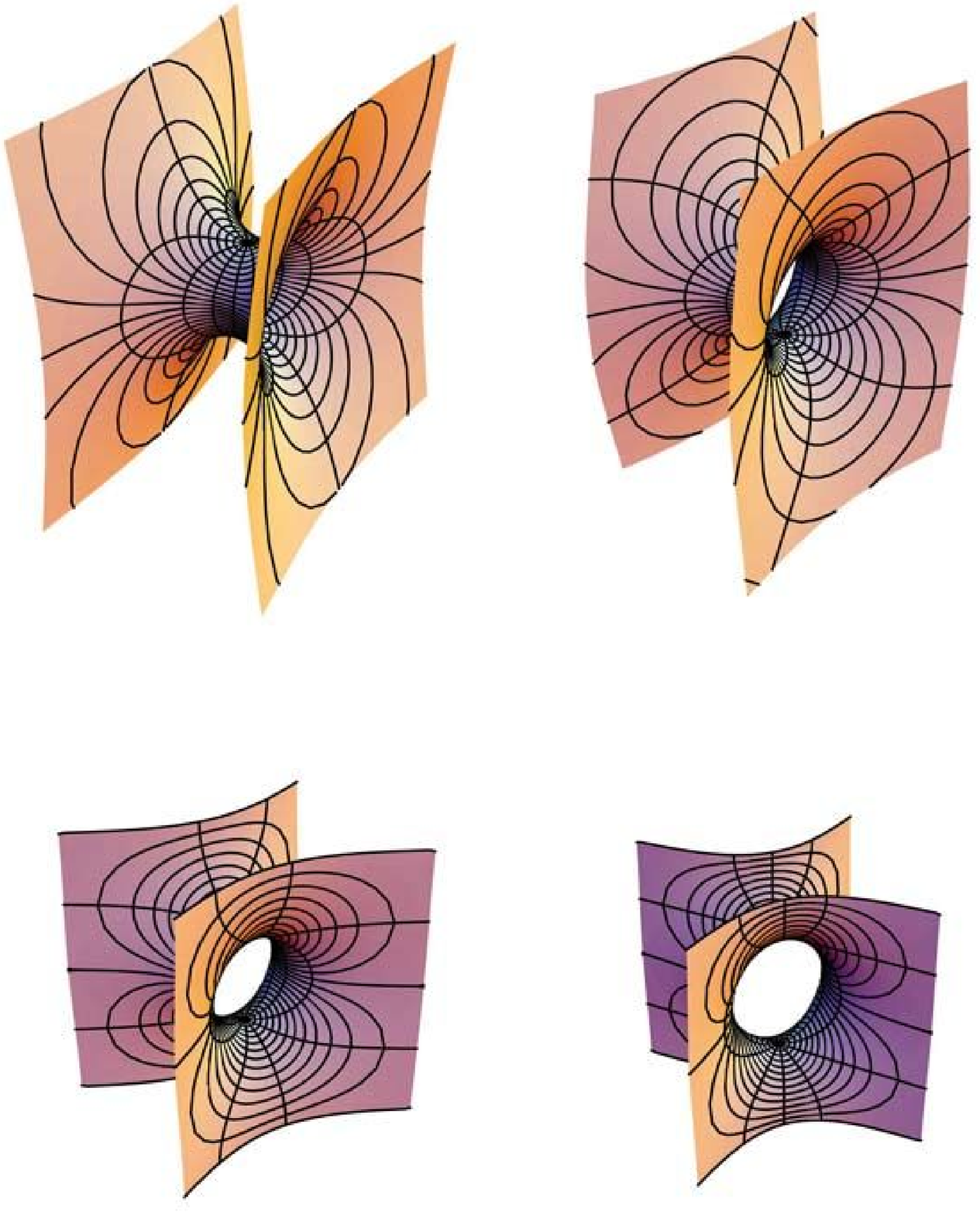}}}
\caption{Ennepers to Scherks singly-periodic transformation for
$t=i/5$ for $i=0,\ldots 5$.} \label{enntosch}
\end{figure}
\end{center}
\end{ex}

%
%
%
%
%
%
%
%
%
%
%
%

\section{Linear combinations that are not convex in on direction}

\begin{ex}[The 4-noid to 4-Enneper]

The harmonic function that lifts to the 4-ennepers surfaces is given by

$$h_{4E}+g_{4E}=z-\frac{z^7}{7}\qquad h_{4E}-g_{4E}=z+\frac{z^7}{7}$$

and that of the 4-noid is given by

$$h_{4N}+g_{4N}=\frac{1}{8}\left(\frac{2z}{1+z^2}-3\log\left(\frac{z+1}{z-1}\right)\right) \qquad h_{4N}-g_{4N}=
\frac{1}{4}\left(\frac{z}{1-z^2}+\frac{3i}{2}\log\left(\frac{1-iz}{1+iz}\right)\right)$$

We restrict the domain of these surfaces to $B(0,.95)$ to avoid self intersections.
Neither of these surfaces is convex, thus we need to pursue other means then the above for showing that the combination
$f=s f_{4E}+(1-s)f_{4N}$ is minimal for all $s\in(0,1)$.  The following lemma will prove univalence: \end{ex}

\begin{lem} Let $s$ be fixed such that $0 \leq s < 1$. For any $n \geq 2$, $f$ is univalent in $\D$. \end{lem}

\begin{proof} Fix $r_0$ such that $0<r_0<1$ and consider $\Omega \subset \D$ the region bounded by $\sigma_1 \cup
\{0\}, \sigma_2, \sigma_3,$ and $\sigma_4$, where $\sigma_1=\{r:0 < r \leq r_0 \}$, $\sigma_2=\{re^{i \pi /4}:0 < r
\leq 1 \}$, $\sigma_3= \{e^{i \pi (1-r)/4}:0 \leq r \leq r_0\}$, and $\sigma_4=\{ z=tr_0+(1-t)e^{i \pi (1-r_0)/4}:0
\leq t \leq 1\}$. We will prove this claim in three steps. First, we will show that $f$ is univalent in $\Omega$ for
$r_0$ arbitrarily close to 1, and that $0 \leq Arg(f(\Omega)) \leq \frac{\pi}{4}$. Second, we verify that $f$ is
univalent in the sector $\Omega \cup \Omega'$, where $\Omega'$ is the reflection of $\Omega$ across the real axis, and
$\frac{- \pi}{4} \leq Arg(f(\Omega \cup \Omega')) \leq \frac{\pi}{4}$. Finally, we will verify that $f$ is univalent in
$\D$.

Step One: The argument principle for harmonic functions \cite {durenarg} is valid if $f$ is continuous on
$\overline{D}$, $f(z) \neq 0$ on $\partial D$, and $f$ has no singular zeros in $D$, where $D$ is a Jordan domain. Note
$z_0$ is a singular point if $f$ is neither sense-preserving nor sense-reversing at $z_0$. We will show that for arbitrary $M > 0$, we may choose $r_0 <
1$ so that each value in the region bounded by $|w| < M$ and $0 < Arg(w) < \frac{\pi}{4}$ is assumed exactly once in
the sector bounded by $|z| < 1$ and $0 < Arg(z) < \frac{\pi}{4}$, while no value in the region bounded by $|w| < M$ and
$\frac{\pi}{4} < Arg(w) < 2 \pi$ is assumed in this sector.

Observe that $f_1'(z)=0$ only if $z$ is an 4th root of -1. Thus, on $\sigma_1$, $f_1$ is an increasing function of $r$
with $Arg(f_1) =0$. Also, as $|z|$ increases on $\sigma_2$ and $Arg(z)$ decreases on $\sigma_3$, $|f_1(z)|$ increases.
Note that $Arg(f_1(\sigma_2 \cup \sigma_3))= \frac{\pi}{4}$.  For $f_2$, if we let $z=\rho e^{i \theta}$ and use the
fact that $f_2=h_2+\overline{g_2}$, we get

Note that $f_2(0)=0$. For $z \in \sigma_1$, $\frac{d}{d \rho} (f_2(\rho)) > 0$, and so $f_2$ increases on $\sigma_1$ as $r$ increases. Also $f_2(\rho) > 0$; hence
$Arg(f_2(\sigma_1))=0$. For $z \in \sigma_2$, $\frac{d}{d
\rho} (f_2(\rho e^{i \pi /4})) \neq 0$, and so $f_2
(\sigma_2)$ does not reverse its direction. Finally we note $Arg(f_2(\sigma_2))=\frac{\pi}{4}$. Recall $f_2$ is constant on $\sigma_3$.  Therefore, we
see that for $j=1,2,3$, $f(\sigma_j)$ is a simple curve with $Arg(f (\sigma_1)) =0$ while $Arg(f(\sigma_2 \cup
\sigma_3))=\frac{\pi}{4}$. To complete the proof that $f$ is univalent on $\Omega$, it suffices to show that given any
$M>0$ there exists an $r_0$ such that $|f(z)|>M$ for all $z \in \sigma_4$. To see this note that $|f_2(z)|$
is bounded  for all $z \in \D$ while for $s$ fixed $(0 \leq s < 1)$ and for $z \in \sigma_4$,~~$(1-s)f_1(z) \rightarrow
\infty$ as $r \rightarrow 1$. Hence for a given $M$ the inequality will hold if we take $r_0$ sufficiently close to
$1$. The proof is now complete since we have shown that every point outside the wedge is not assumed while every point
inside the wedge is assumed exactly once by $f$.

Step Two: Since $f$ is univalent in $\Omega$, we can use reflection across the real axis to establish that $f$ is
univalent in the sector $\Omega '$. In particular, suppose $z_1,z_2 \in \Omega '$ with $f(z_1)=f(z_2)$. Then by
symmetry $\overline{f(\overline{z_1})} =f(z_1)=f(z_2)=\overline{f(\overline{z_2})}$. Hence, $f(\overline
{z_1})=f(\overline{z_2})$, or $\overline{z_1}=\overline{z_2}$. Arguing in the same manner as in Step One, we can show
that $0 \geq Arg(f(\Omega ')) \geq \frac{- \pi}{4}$. Therefore, $f$ is univalent in $\Omega \cup \Omega '$ and its
image is in the wedge between the angles $\frac{- \pi}{4}$ and $\frac{\pi}{4}$.

Step Three: First, it is true that $e^{i \pi j/2} f(ze^{-i \pi j/2}) =f(z)$, for all $z \in \D$ where $j=0,1,...,4$.
To see this note that

Now, using this fact that $e^{i \pi j/2} f(ze^{-i \pi j/2})=f(z)$, we see that if $z$ is any point in
$\D$, it can be rotated so that it is in the sector $\Omega '$, in which $f$ is univalent, and then rotated back by
multiplying by the constant $e^{i \pi j/2}$ and hence preserving univalency. \end{proof}

\begin{ex}[The 4-noid to 4-Enneper]

The harmonic function that lifts to the 4-ennepers surfaces is given by

$$h_{4E}+g_{4E}=z-\frac{z^7}{7}\qquad h_{4E}-g_{4E}=z+\frac{z^7}{7}$$

and that of the 4-noid is given by

$$h_{4N}+g_{4N}=\frac{1}{8}\left(\frac{2z}{1+z^2}-3\log\left(\frac{z+1}{z-1}\right)\right) \qquad h_{4N}-g_{4N}=
\frac{1}{4}\left(\frac{z}{1-z^2}+\frac{3i}{2}\log\left(\frac{1-iz}{1+iz}\right)\right)$$

Neither of these surfaces is convex, thus we need to pursue other means then the above for showing that the combination
$f=s f_{4E}+(1-s)f_{4N}$ is minimal for all $s\in(0,1)$. \end{ex}

\begin{center}
\begin{figure}[h]
\centerline{\hbox{
\includegraphics[height=5in]{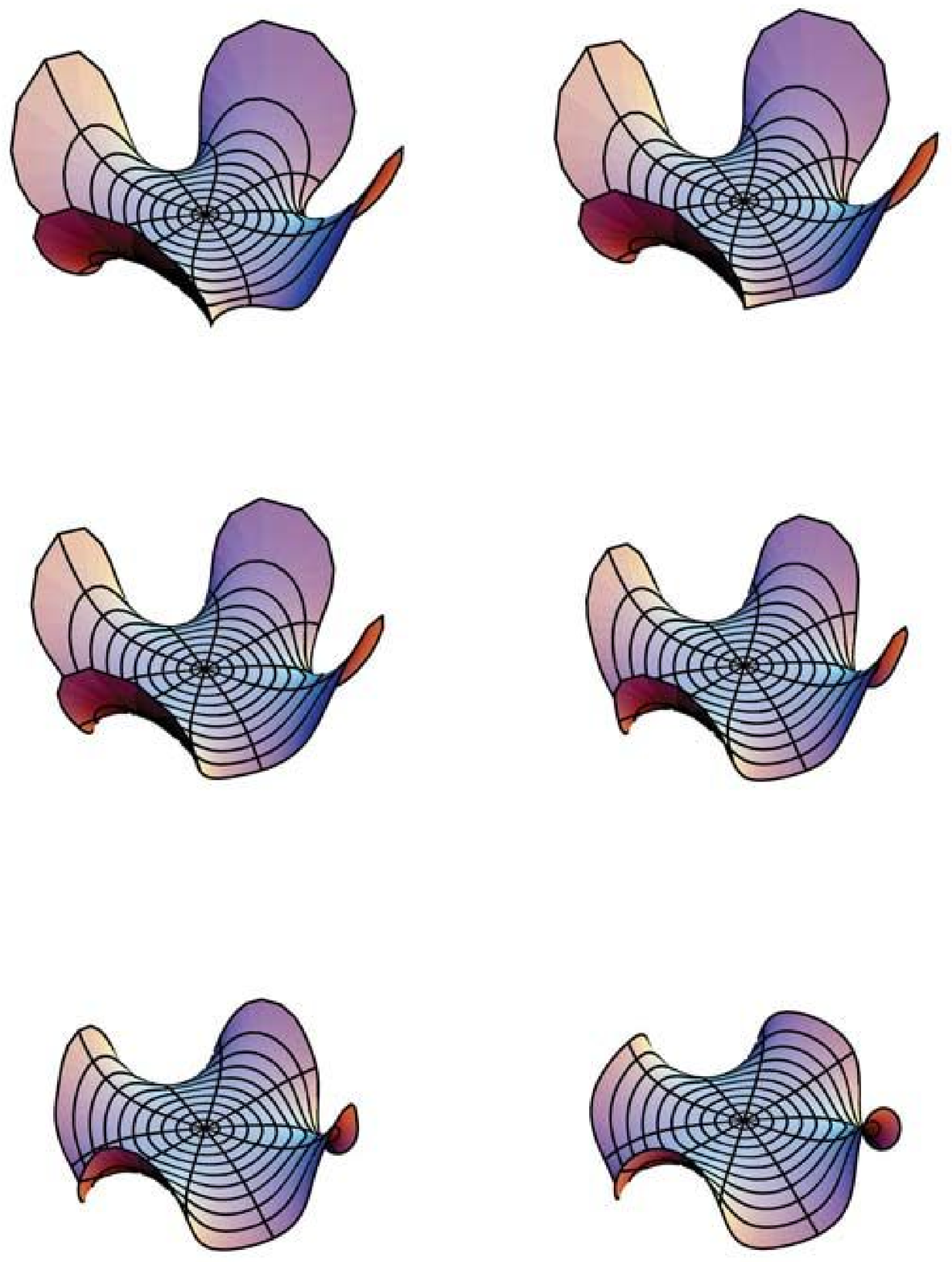}}}
\caption{\label{fig2}  4 Enneper to 4-noid}
\end{figure}
\end{center}

We plot in $(\ref{fig2})$ four instances of this transformation.

\begin{thebibliography}{10}

\bibitem{bdp} R. Berry, M. Dorff, and W. L. Petersen, Lie
symetries of minimal surfaces, preprint.

\bibitem{css}
J. Clunie and T. Sheil-Small, Harmonic univalent functions, {\em
Ann. Acad. Sci. Fenn. Ser. A.I 074 Math.} {\bf 9} (1984), 3-25.

\bibitem{dhkw} U. Dierkes, S. Hildebrandt, A. Küster, and O. Wohlrab,
{\it Minimal surfaces I}, Grundlehren der Mathematischen
Wissenschaften, 295, Springer-Verlag, Berlin, 1992.

\bibitem{dorffszynal}
M. Dorff and J. Szynal, Harmonic shears of elliptic integrals,
{\em Rocky Mountain Journal of Mathematics}, {\bf 35} (2005), no.
2, 485-499.

\bibitem{driverduren} K. Driver and P. Duren,
Harmonic shears of regular polygons by hypergeometric functions,
{\em J. Math. Anal. App.} {\bf 239} (1999), 72-84.

\bibitem{Duren}
P. Duren, {\it Harmonic mappings in the plane}, Cambridge Tracts
in Mathematics, 156, Cambridge University Press, Cambridge, 2004.

\bibitem{dms}
P. Duren, J. McDougall, and L. Schaubroeck, Harmonic mappings onto
stars, {\it J. Math. Anal. Appl.} 307 (2005), no. 1, 312-320.

\bibitem{hs}
W. Hengartner and G. Schober, On schlicht mappings to domains
convex in one direction,  {\it Comment. Math. Helv.} 45 (1970),
303-314.

\bibitem{convex}
W. Koepf, Parallel accessible domains and domains that are convex
in some direction, {\it Partial Differential Equations with
Complex Analysis}, Pitman Res. Notes Math. Ser., 262, Longman Sci.
Tech., Harlow, 1992, 93-105.

\bibitem{lew} H. Lewy,  On the non-vanishing of the Jacobian in
certain one-to-one mappings, {\it Bull. Amer. Math. Soc.}
{\bf 42} (1936), 689-692.

\bibitem{convextwo}
W. C. Royster and M. Ziegler, Univalent functions convex in one
direction, {\it Publ. Math. Debrecen.}, \bf{23}, no. 3-4, 339-345.


\end{thebibliography}
\end{document}